\documentclass[a4paper,10pt]{article}

\usepackage{amsfonts}
\usepackage{amsmath}
\newtheorem{thm}{Theorem}
\newtheorem{lem}[thm]{Lemma}
\newtheorem{cor}[thm]{Corollary}
\newtheorem{dfn}[thm]{Definition}
\begin{document}
\title {Spectrum of the Laplacian on radial graphs}
\author{Rodrigo Bezerra de Matos and Jos\'e Fabio B. Montenegro\thanks{\small rodrigolmatos93@gmail.com, fabio@mat.ufc.br }}
\date{\today}
\maketitle
\begin{abstract}We prove that if $M$ is a complete hypersurface in $\mathbb{R}^{n+1}$ which is graph of a real radial function, then
the spectrum of the Laplace operator on M is the interval $[0,\infty)$.

\vspace{.2cm}

\vspace{.2cm}
\end{abstract}

\section{Introduction}
Let $M$ be a simply connected Riemannian manifold. The Laplace operator $\Delta :C_{0}^{\infty}(M)\to C_{0}^{\infty}(M)$, defined as
$\Delta =\mathrm{div}\circ\mathrm{grad}$ and acting on $C_{0}^{\infty}(M)$, the space of smooth functions with compact support,
 is a second order elliptic operator and it has a unique extension $\Delta$ to an unbounded self-adjoint operator
on $L^2(M)$, whose domain is $\mathrm{Dom}(\Delta)=\{f\in L^{2}(M): \Delta f \in L^{2}(M)\}$. Since $-\Delta$ is positive and symmetric, its spectrum is the set of $\lambda \geq0$ such that
$\Delta +\lambda I$ does not have bounded inverse. Sometimes we say spectrum of M rather than spectrum of $-\Delta$  and we denote it by $\sigma(M)$.
One defines the {\it essential spectrum} $\sigma_{ess}(M)$ to be
those $\lambda$ in the spectrum which are either accumulation points of the spectrum or eigenvalues of infinite multiplicity.
The {\it discrete spectrum} is the set $\sigma_d=\sigma (M)\setminus \sigma_{ess}(M)$ of all eigenvalues of finite multiplicity which are isolated point of the spectrum.

There is a vast literature studying  the spectrum of the Laplace operator on complete non-compact manifolds.
See \cite{bessa-Jorge-montenegro-JGA}, \cite{bessa-Jorge-mari}, \cite{bmp}, \cite{DoLi}, \cite{kleine1}, \cite{kleine2},  for  geometric conditions implying the discreteness of the spectrum,
$\sigma_{\mathrm{ess}}(M)=\emptyset$. For the discreteness of the spectrum of bounded minimal submanifolds, one can
see \cite{bessa-Jorge-montenegro-JGA}, \cite{bessa-Jorge-mari}. For purely continuous spectrum ($\sigma_{d}(M)=\emptyset)$ see \cite{Chen},
\cite{donnelly}, \cite{donnelly-garofalo}, \cite{Escobar},  \cite{karp}, \cite{Kumura}, \cite{Li}, \cite{Zhiqin and Zhou}, \cite{tayoshi}, \cite{wang}, \cite{Zhou}.
In this work we consider geodesically complete hypersurfaces
which are graphs of a radial functions. Our main result is the following theorem.

\begin{thm}
\label{thm1}
Let $M$ be a complete hypersurface in $\mathbb{R}^{n+1}$ which is graph of a real radial function.
Then, the spectrum of the Laplace operator on M is $[0,\infty)$.
\end{thm}

The Theorem above allow us to construct a bounded hypersurface with the same spectrum of $\mathbb{R}^{n+1}$.
For instance, by taking $M$ to be the graph of the function $f(x)=\cos( \tan (\pi |x|/2))$ defined on the unit open ball.
The end of this hypersurface is somewhat chaotic, in terms of convergence of the curvatures or mean curvature of geodesics spheres.
As a consequence, one cannot use the existing results to determine its spectrum.

Since M is complete, and  graph of an radial function $f:D\to \mathbb{R}$,
D is an open ball or $D=\mathbb{R}^{n}$. We consider on D the spherical coordinate system $X:[0,R)\times{\cal{O}}\to D$,
defined by $X(r,x_1,...,x_{n-1})=r\,w(x_1,...,x_{n-1})$, where $0<R\leq +\infty$ and $w$ is a coordinate system on $S^{n-1}$ defined on an open set $\cal{O}$ of $\mathbb{R}^n$.
Note that M has a natural coordinate system $Y:[0,R)\times\cal{O}\to M$, given by $Y(r,x_1,...,x_{n-1})=( r\,w(x_1,...,x_{n-1}),f(r))$, but we are interested on
the spherical coordinate system for M on $p=(0,f(0))$.
For this, consider  $t:[0,R)\to [0,\infty)$ given by$$t(r)=\int_{0}^{r}\!\sqrt{1+f'(\tau)^2}\,d\tau$$
and  observe that $$t'(r)=\sqrt{1+f'(r)^2}>0$$ with
$$\lim_{r \to R}t(r)=+\infty$$

so that $t:[0,R)\to[0,\infty)$ is an diffeomorphism. We denote by $r:[0,\infty)\to[0,R)$ the inverse diffeomorphism.
By the inverse function theorem,
$$0<r'(t)=\frac{1}{\sqrt{1+f'(r(t))^2}}\leq1.$$

Finally, the system of spherical coordinates on M,  $Z:[0,\infty)\times\cal{O}\to M$ is defined by
$$Z(t,x_1,...,x_{n-1})=(r(t)\,w(x_1,...,x_{n-1}),f\!\circ\! r(t)).$$
The metric of M on such system is given by  $$g_M=dt^2+r(t)^2g_{S^{n-1}}$$

This way, Theorem \ref{thm1} is a rather simple consequence of theorem below about a class of model manifolds.

\begin{thm}
\label{thm2}
Let $M=[0,\,\infty )\times \mathbb{S}^{n-1}$ be a model manifold with metric given by $g_M=dt^2+r^2(t)g_{\mathbb{S}^{n-1}}$ such that $0<r'(t)\leq c$ for all $t\geq 0$.
Then, the spectrum of the Laplace operator on M is $[0,\infty)$.
\end{thm}
On the next section we prove the Theorem 2. The last section we deal with the Sturm-Liouville Theory.

\section{Proof of Theorem 2}
Since $r'(t)>0$, $r(t)$ is increasing and there exists only two possibility $$\lim_{t\to \infty}r(t)=\infty \;\;\mathrm{or} \;\;\lim_{t\to \infty}r(t)=R.$$
 In the first case, since $r'(t)$ is limited, we have
$$\lim_{t\to\infty}\Delta\, t=\lim_{t\to\infty}\frac{r'(t)}{r(t)}=0.$$
By the result \cite[Theorem 1.2]{Kumura}, obtained by H. Kumura, it follows that the spectrum of $M$ is the interval $[0,\infty )$.

Let us now turn to the case where $\lim_{t\to \infty}r(t)=R.$ To start, we present a well know characterization for the spectrum
of a self-adjoint operators in a Hilbert space \cite[Lemma 4.1.2, p.73]{Davies}
\begin{lem}\label{lem3}
A number $\lambda \in \mathbb{R}$ lies in the spectrum of a self-adjoint operator $H$ if and only if there exists a sequence of functions
$f_n\in \mathrm{Dom}H$ with $\|f_n\|=1$ such that
$$ \lim_{n\to \infty}\|Hf_n-\lambda f_n\|=0.$$

\end{lem}

To deduce theorem \ref{thm2} from lemma \ref{lem3} we will construct, for each $\lambda > $0, a sequence of radial smooth functions $f_p:M\to \mathbb{R}$ with compact support such that
\begin{equation}
\label{desig-f}
\|\Delta f_p+\lambda f_p\|_{L^ 2(M)}\leq \frac{c}{p}\, \|f_p\|_{L^ 2(M)}
\end{equation}
for any natural $p$, where $c$ is a constant which does not depend on $p$. So that $g_p=f_p/\|f\|_2$ has norm one and
$$\lim_{p\to \infty}\|\Delta g_p+\lambda g_p\|_{L^ 2(M)}=0.$$Therefore, by Lemma \ref{lem3}, $\lambda $ lies spectrum.
To construct the function $f_p$, we fix $t_0> 0$ and prove that there are $t_{1}(\lambda)> t_0$ and a radial function $u = u (t)$ solution of the problem
\begin{equation}
\label{laplace-u}
\left\{
\begin{array}{rll}
\Delta u+\lambda\,u&=&0 \;\;\;\mathrm{in}\; [t_0,\,t_1] \\
u(t_0)=u(t_1)&=&0 \\
u&>&0 \;\;\;\mathrm{in}\; (t_0,\,t_1)
\end{array}
\right.
\end{equation}
Using the Sturm-Liouville theory, we showed that $u$  can be extended to the whole interval $[t_0, \infty )$ and it has infinite zeros $t_0 <t_1<\dots < t_p<\dots$.
The next step is to consider, for each $p$, a smooth bump function $h_p$ whose support is the interval $[t_0, \, t_{3p}]$,
define $f_p=u\,h_p$ and show that each $f_p$ in this sequence satisfies (\ref{desig-f}).

First we observe that the equation in (\ref{laplace-u}) is equivalent to
 \begin{equation}
  \label{equ}
 [r^{n-1}(t)u'(t)]'+\lambda\,r^{n-1}(t)u(t)=0
 \end{equation}
if $u=u(t)$ is a radial function. By Theorem \ref{thm12}, see the Appendix, for every $\lambda>0$ there is a function $u\colon[t_0,\infty)\mapsto\mathbb{R}$ solution of
(\ref{equ}) with $u(t_0)=0$.

By \cite[Thm. 2.1]{Manfredo}, see also Corollary \ref{cor6},  $u$ has arbitrarily large zeros. Let $t_0<t_1<\dots$  the increasing sequence of zeros for $u$.

 For $p \in \mathbb{N}$, we choose a smooth bump function $h=h_p\colon\mathbb{R}\mapsto \mathbb{R}$ such that
 $$\left\{
\begin{array}{l}
 h(t)=0 , \, t\in(-\infty,t_0]\cup [t_{3p},\infty) \\ \\
 h(t)=1, \, t\in [t_{p},{t_{2p}}]\\ \\
0\leq {h} \leq1
\end{array}\right.$$
For instance, the function $h$ is defined in the following way: let $\varphi\in C_{0}^{\infty}(\mathbb{R})$ be a function satisfying
 $ \varphi \geq 0$, $ \mathrm{supp}{\varphi}=[0,1]$ and $\int\varphi =1$.
Now, put $$h_{p}(t)=\int_{-\infty}^{t} \varphi_{p}(s)ds$$
where $$\varphi_{p}(t)=\frac{1}{t_p-t_0}\varphi \left (\frac{t-t_0}{t_p-t_0}\right)- \frac{1}{t_{3p}-t_{2p}}\varphi \left(\frac{t-t_{2p}}{t_{3p}-t_{2p}} \right).$$
We observe that $h$  satisfies
\begin{equation}
 \label{esth}
 \begin{array}{l}
 \| h_{p}'\| _{\infty}\leq \mathrm{max}\left\{{\displaystyle\frac{\|\varphi\|_{\infty}}{t_p-t_0}, \frac{\|\varphi\|_{\infty}}{t_{3p}-t_{2p}} }\right\} \leq \displaystyle\frac{c}{p} \\ \\
 \| h_{p}''\| _{\infty}\leq \mathrm{max}\left\{{\displaystyle\frac{\|\varphi'\|_{\infty}}{(t_p-t_0)^2}, \frac{\|\varphi'\|_{\infty}}{(t_{3p}-t_{2p})^2} }\right\} \leq \displaystyle\frac{c}{p^2}
 \end{array}
\end{equation}
since $(t_{i+1}-t_i)\geq C$, for all $i\geq 0$, and  some universal constant $C$, according to Corollary \ref{estzeros} in Section \ref{Apendice}.

Consider $f=f_p=u\cdot h_{p}$ and we  prove that such function satisfies the inequality in (\ref{desig-f}). Computing $\Delta f + \lambda f$ we obtain,

 $$\Delta f+ \lambda f=2u'h'+uh''+(n-1)\frac{r'}{r}h'u$$
 Using the inequalities (\ref{esth}),  $r(t)\geq  r(t_0)$, $0< r' \leq c$ we have

 $$|\Delta f+ \lambda f|\leq \frac{c}{p}(|u'|+|u|)\chi_{[t_0,t_{3p}]}$$
 Then,
 $$|\Delta f+ \lambda f|^2\leq \frac{c}{p^2}(|u'|^2+|u|^2)\chi_{[t_0,t_{3p}]}$$

 $$ \int_M|\Delta f+ \lambda f|^2dM\leq \frac{c}{p^2}\left(\int_{t_0}^{t_{3p}}|u'|^2r^{n-1}dt+\int_{t_0}^{t_{3p}}|u|^2r^{n-1}dt\right)$$

Mutiplying (\ref{equ}) by $u$, integrating, and using the integration by parts formula we find
$$\int_{t_0}^{t_{3p}}|u'|^2r^{n-1}dt=\lambda\int_{t_0}^{t_{3p}}|u|^2r^{n-1}dt$$

$$\|\Delta f_p+ \lambda f_p\|_{L^2(M)}\leq \frac{c}{p}\;\| u.\chi_{[t_0,t_{3p}]}\|_{L^2(M)}
\leq \frac{c}{p}\;\| u.\chi_{[t_p,t_{2p}]}\|_{L^2(M)}\leq \frac{c}{p}\;\|f_p\|_{L^2(M)}$$
where the second inequality comes from Lemma \ref{lemma12} below.

\begin{lem}
\label{lemma12}
There is a positive constant $c$ independent on $p$ such that
$$\int_{t_0}^{t_{3p}}u^2r^{n-1}dt\leq c\,\int_{t_p}^{t_{2p}}u^2r^{n-1}dt$$
where $u$ is solution of (\ref{equ}) and $t_0<t_1<\dots$ are zeros of $u$.
\end{lem}

Observe that, multiplying  (\ref{equ}) by $r^{n-1}u'$ we get $$ (r^{n-1}u')'r^{n-1}u'+\lambda r^{2(n-1)}uu'=0$$ and so, $$ [(r^{n-1}u')^2]'+\lambda r^{2(n-1)}(u^{2})'=0 $$
Integrating from $t_0$ to $t_k$ we have
$$r(t_k)^{2(n-1)}u'(t_k)^{2}-r(t_0)^{2(n-1)}u'(t_0)^2=-\lambda\int_{t_0}^{t_k}\!\!r^{2(n-1)}(s)[u^{2}(s)]'ds$$
Since the right side is equal to$$2\lambda(n-1)\!\int_{t_0}^{t_k}\!r^{2(n-1)-1}r'u^{2}ds$$
we find
\begin{equation}
\label{eq1}
r(t_k)^{2(n-1)}u'(t_k)^{2}-r(t_0)^{2(n-1)}u'(t_0)^2=2\lambda(n-1)\!\int_{t_0}^{t_k}\!r^{2(n-1)-1}r'u^{2}ds
\end{equation}

  Provided the right side of the last equation is positive and $r<R$,
  \begin{equation}
  \label{est1}
    u'(t_k)^2>\frac{r(t_0)^{2(n-1)}}{R^{2(n-1)}}u'(t_0)^2
    \end{equation}
    for any $k\geq1$.

 Now to obtain an estimative in the other direction we observe that for
  $$w_0(t)= u'(t_0)\lambda^{-1/2}\sin\left(\sqrt{\lambda} \;r^{n-1}(t_0)\int_{t_0}^{t}\frac{ds}{r^{n-1}(s)}\right)$$ we have $w'(t_0)=u'(t_0)>0$ and
 \begin{equation}
  \label{eq2}
(r^{n-1}(t)w'(t))'+\frac{\lambda r(t_0)^{2(n-1)} w(t)}{r^{n-1}(t)}=0
 \end{equation}
 Multiplying by $r^{n-1}(t)w'$ we get, similarly to the calculations above
 \begin{equation}
 \label{eqw}
  [r(t)^{2(n-1)} (w')^2]'+ \lambda r(t_0)^{2(n-1)}(w^{2})'=0
  \end{equation}
 Now, if $\overline{t_1}$ is the next root of $w$ after $t_0$, integrating the last equation we find
 \begin{equation}
 \label{eq3}
r(\overline{t_1})^{2(n-1)}w'(\overline{t_1})^2=r(t_0)^{2(n-1)} w'(t_0)^2=r(t_0)^{2(n-1)} u'(t_0)^2
 \end{equation}
 Now, the right side of (\ref{eq1}) , for $k=1$ can be estimated in the following way $$ 2(n-1)\lambda\int_{t_0}^{t_1}r^{2n-3}r'u^2\leq2(n-1)\lambda\int_{t_0}^{t_1}r^{2n-3}r'w^2
 \leq2(n-1)\lambda\int_{t_0}^{\overline{t_1}}r^{2n-3}r'w^2$$

 $$=\lambda\int_{t_0}^{\overline{t_1}}(r^{2(n-1)})'w^2
 =-\lambda\int_{t_0}^{\overline{t_1}}r^{2(n-1)}(w^2)'$$

 $$=-\frac{1}{r(t_0)^{2(n-1)}}\int_{t_0}^{\overline{t_1}}r^{2(n-1)}(\lambda r(t_0)^{2(n-1)} w^2)'$$
 And by (\ref{eqw}) the last term is equal to
 $$\frac{1}{r(t_0)^{2(n-1)}}\int_{t_0}^{\overline{t_1}}r^{2(n-1)}(r^{2(n-1)} (w')^2)'$$
$$ =\frac{1}{r(t_0)^{2(n-1)}}\int_{t_0}^{\overline{t_1}}(r^{4(n-1)}(w_t)^2)' - 2(n-1)r^{4n-5}r'(w_t)^{2}$$

 $$<\frac{r^{4(n-1)}(\overline{t_1})(w_t)^{2}(\overline{t_1})-r^{4(n-1)}(t_0)(w_t)^{2}(t_0)}{r(t_0)^{2(n-1)}}$$

 Now, using (\ref{eq3}) and that $w'(t_0)=u'(t_0)$, we find $$2(n-1)\lambda\int_{t_0}^{t_1}r^{2n-3}r'u^2\leq(r(\overline{t_1})^{2(n-1)}-r(t_0)^{2(n-1)})u'(t_0)^{2}$$

 Then, by (\ref{eq1}),
 $$r(t_1)^{2(n-1)}u'(t_1)^{2}-r(t_0)^{2(n-1)}u'(t_0)^{2}\leq(r(\overline{t_1})^{2(n-1)}-r(t_0)^{2(n-1)})u'(t_0)^2 $$
 Since $r(t)$ is increasing, it follows that
 $$ r(t_1)^{2(n-1)}u'(t_1)^2\leq r(\overline{t_1})^{2(n-1)}u'(t_0)^2\leq r(t_2)^{2(n-1)}u'(t_0)^2$$

 Then, $$u'(t_1)^{2}\leq \frac{r(t_2)^{2(n-1)}}{r(t_0)^{2(n-1)}} u'(t_0)^2$$ Using the same argument, one shows by induction that
 $$u'(t_k)^{2}\leq \frac{r(t_{k+1})^{2(n-1)} r(t_k)^{2(n-1)}}{r(t_1)^{2(n-1)} r(t_0)^{2(n-1)}} u'(t_0)^2$$And since $r(t)<R$ we find that
 \begin{equation}
 \label{estimativa}
 u'(t_k)^{2}\leq \frac{R^{4(n-1)}}{r(t_0)^{2(n-1)} r(t_1)^{2(n-1)}}u'(t_0)^{2}
\end{equation}
 Now, using corollary \ref{wv} its easy to check that

 $$\int_{t_0}^{t_{3p}}u^2r^{n-1}=\sum_{k=0}^{3p-1}\int_{t_k}^{t_{k+1}}u^2r^{n-1}(t)dt$$
 $$\leq \frac{1}{\lambda}\sum_{k=0}^{3p-1}u'(t_k)^2\int_{t_k}^{t_{k+1}}\sin^2 \left(\sqrt{\lambda}\;r^{n-1}(t_k)\int_{t_k}^{t}\frac{ds}{r^{n-1}(s)} \right)r^{n-1}(t)dt$$
 Consequently, if $$\tau=\sqrt{\lambda}\;r^{n-1}(t_k)\int_{t_k}^{t}\frac{ds}{r^{n-1}(s)}$$ the change of variables formula shows that the last sum is equal to
  $$\frac{1}{\lambda^{3/2}}\sum_{k=0}^{3p-1}\frac{u'(t_k)^2}{r^{n-1}(t_k)}
  \int_{0}^{\pi}\sin^2(\tau)r^{2(n-1)}(\tau(t))d\tau$$
 $$  \leq \frac{\pi R^{2(n-1)}}{2\lambda^{3/2}r^{n-1}(t_{0})}\sum_{k=0}^{3p-1}u'(t_k)^2=C\sum_{k=0}^{3p-1}u'(t_k)^2$$
  Since by the equations (\ref{est1}) and (\ref{estimativa}) the following inequalities hold
  $$\sum_{k=0}^{3p-1}u'(t_k)^2\leq 3C\,p u'(t_0)^{2}
 \leq C\sum_{k=p}^{2p-1}u'(t_k)^2$$
  we have
 \begin{equation}
 \label{compsom}
 \int_{t_0}^{t_{3p}}u^2r^{n-1}\leq C\sum_{k=p}^{2p-1}u'(t_k)^2
 \end{equation}
 Where the last inequality comes from (\ref{est1}), for some suitable constant $C>0$.
 Again by the change of variables formula, this time applied to each $v_k$ and by corollary \ref{wv} one sees that if
 $\tilde t_k$ is the next zero of $v_k$ after $t_k$ we have
  $$\int_{t_p}^{t_{2p}}u^2r^{n-1}(t)dt=
  \sum_{k=p}^{2p-1}\int_{t_k}^{ t_{k+1}}u^2r^{n-1}(t)dt$$
  $$\geq \sum_{k=p}^{2p-1} \int_{t_k}^{\tilde t_{k+1}}v_{k}^2\,r^{n-1}(t)dt$$
  $$\geq{C}\sum_{k=p}^{2p-1}u'(t_k)^2$$
  From (\ref{compsom}) we conclude that
  $$\int_{t_0}^{t_{3p}}u^2r^{n-1} \leq C\int_{t_p}^{t_{2p}}u^2r^{n-1}$$
for every $p\in\mathbb{N}$ and for a constant ${C}={C}(\lambda,R)$ independent of $p$.

\section{Appendix: \large{Elements of Sturm-Liouville Theory}}\label{Apendice}

Since for a radial function $u$
$$\Delta u +\lambda u=\frac{(r^{n-1}(t)u')'+ \lambda r^{n-1}(t)u}{r^{n-1}(t)}$$
 in order to study eigenvalues and the spectrum of $\Delta$ is natural to ask some information about the Sturm-Liouville problem

\begin{equation}
 \label{sturm}
(r^{n-1}(t)u')'+\lambda r^{n-1}(t)u=0 , \;\;t\geq t_0>0
\end{equation}
Or, more generally, about

\begin{equation}
\label{sturm2}
(v(t)u')'+\lambda v(t)u=0 , \;\;t\geq t_0>0
\end{equation}
where $v(t)$ a positive continuous function on $[t_0,\infty)$ and $\lambda > 0$.

We start our study with a classical terminology

\begin{dfn}The equation (\ref{sturm2}) is said to be oscillatory if any of it's solutions  has arbitrarily large zeros.
\end{dfn}

The following theorem is a practical criterion for oscillation.
\begin{thm}
Let $v(t)$ a positive continuous function on $[t_0,\infty)$ and $\lambda > 0$.
Then, the equation $$(v(t)u')'+\lambda v(t)u=0$$ for $t\geq t_0$
is oscillatory, provided $\int_{t_0}^{\infty}v(t)dt=+\infty$ and
$\int_{t_0}^{t}v(t)dt\leq C\;t^a$, for some positive constants C and a.

\end{thm}
The proof is discussed on \cite[Theorem 2.1]{Manfredo}.
\begin{cor}
\label{cor6}
The equation (\ref{sturm}) is oscillatory.
\end{cor}
\textbf{Proof}:
Since $$\lim_{t\to\infty}\! r(t)=R$$
 there is a $t_1>0$ such that
$R/2<r(t)<R$ for $t>t_1$. Consequently, $$(R/2)^{n-1}<r^{n-1}(t)<R^{n-1}$$
 for $t>t_1$ and so the criterion above is applied, immediately.

The next theorems emerge as a useful way to compare solutions for different ordinary equations and will be used in the next paragraph.
\begin{thm}
\label{thm8}
Let $x,\,y$ non-trivial solutions for
$$\left\{
\begin{array}{l}
 (p(t)x')'+ q(t)x=0 \\ \\
 (p_1(t)y')'+q_1(t)y=0
\end{array}\right.$$
 Where $p(t)\geq p_1(t)>0$, $q_1(t)\geq q(t)$ for every $t\in I$.
 If $t_1<t_2$ are consecutive zeros for $x$, then $y$ has a zero on $J=(t_1,t_2)$
 or exist $d\in \mathbb{R}$ for which $y=d\, x$ on $J$.
 \end{thm}
 \textbf{Proof}
 Observe that if $y$ does not have a zero on $J$, then
  $$\left(x\frac{(p(t)x'y-p_{1}(t)xy')}{y}\right)'=(q_{1}-q)x^2 + (p-p_1)(x')^{2}+\frac{p_1(x'y-xy')^2}{y^2}$$
  Integrating from $t_1$ to $t_2$ we have
  $$\int_{t_1}^{t_2}(q_1-q)x^2 dt +
\int_{t_1}^{t_2}(p-p_1)(x')^{2}dt + \int_{t_1}^{t_2} p_1\frac{(x'y-xy')^2}{y^2}dt=0$$ Then, if $y$ is not multiple of $x$, the wroskian $(xy'-x'y)$
  is nonzero on $J$ and we get a contradiction with the last equation.
 \begin{cor}
 \label{estzeros}
 If u solves (\ref{sturm}), let $\{t_{p};p\in \mathbb{N}\}$ a increasing sequence of zeros for u.
  There is a  constant $C>0$, such that
  $$(t_{p+1}-t_{p})>C$$ for any $ p\in \mathbb{N}$
  \end{cor}
  \textbf{Proof}
   Define, for each $p\in \mathbb{N}$ $$\varphi(t)=\sin(\sqrt{2^{n-1}\lambda}(t-t_{p}))$$
  Then, $\varphi$ has a zero on $t=t_p$ and $$(R/2)^{n-1}{\varphi}''+ \lambda R^{n-1} \varphi=0$$ Now, $(R/2)^{n-1}<r^{n-1}(t)<R^{n-1}$ for $t$ sufficiently large, lets say for $t>c_0$. Consequently,
  if $p$ is sufficiently large, we can apply the theorem \ref{thm8} for $u$ and $\varphi$ to conclude that $\varphi$ next zero is on $(t_{p},t_{p+1}).$

   Since the next zero of $\varphi$ after $t_{p}$ is on $t=t_{p} + \frac{\pi}{\sqrt{2^{n-1}\lambda}}$, we have $$t_{p+1}-t_p>\frac{\pi}{\sqrt{2^{n-1}\lambda}}$$ for $t_{p}>c_0$
  from which the corollary follows .
\begin{thm}
\label{teorema9}
Let $x,\,y$ non-trivial solutions for
$$\left\{
\begin{array}{l}
 (p(t)x')'+ q(t)x=0 \\ \\
 (p_1(t)y')'+q_1(t)y=0
\end{array}\right.$$
On a interval $[a,b]$ where $p\geq p_1>0$, $q_1>q$,
$x(a)=0$. Suppose that $c\in(a,b]$ is such that $x(c)\neq 0$,
$y(c)\neq 0$ and $x$ has the same number of zeros as $y$ on $(a,c)$.

Then $$\frac{p(c)x'(c)}{x(c)}\geq \frac{p_1(c)y'(c)}{y(c)}$$

\end{thm}
\textbf{Proof}
 We only deal with the case where $y$ is different from $d.x$, otherwise there is nothing to prove. Let  $a\!=\!a_{0},..., a_{n}$ the zeros for $x$ on $[a,c)$ and
$b_{0},..., b_{n-1}$ the zeros for $y$ on $(a,c)$.
Then, by theorem 6, we have $$a_{i}<b_{i}<a_{i+1}$$ for $i=0,...,n-1$. Consequently,
$y$ has no zero on $(a_{n},c)$. Now, we can use the same idea from the proof of theorem 6 to conclude that
$$\left((px'y-p_{1}xy')\frac{x}{y}\right)'\geq 0$$
on $(a_{n},c)$.
Integrating both sides from $a_{n}$ to $c$ and using that $x(a_{n})=0$ we get $$(px'y-p_{1}xy')(c)\frac{x(c)}{y(c)}\geq 0$$  and since we can always assume that
$x(c)y(c)>0$, we find $$\frac{p(c)x'(c)}{x(c)}\geq \frac{p_{1}y'(c)}{y(c)}$$ as we want.

\begin{cor}
\label{wv}
Let $u$ be a solution for (\ref{sturm}), and choose $t_k$, $t_{k+1}$  to be consecutive zeros for $u$. Define
 $$v_k(t)=\frac{r^{n-1}(t_k)u'(t_k)}{R^{n-1}\sqrt{\lambda}}\sin\left(\sqrt{\lambda}\;R^{n-1}\int_{t_k}^{t}\frac{ds}{r^{n-1}(s)}\right)$$ and
 $$w_k(t)=\frac{u'(t_k)}{\sqrt{\lambda}}\sin\left(\sqrt{\lambda}\;r^{n-1}(t_k)\int_{t_k}^{t}\frac{ds}{r^{n-1}(s)}\right)$$
 then $|v_k|\leq |u|$ on $(t_k,\tilde t_k)$ and $|u|\leq |w_k| $ on $(t_k,t_{k+1})$
 , where  $\tilde t_k$ is the next zero of $v_k$ after $t_k$.
 \end{cor}
\textbf{Proof}
 Observe that $v_k(t_k)=0$, $v_k'(t_k)=u_k'(t_k)$ and
$$(r^{n-1}(t)v'_k)'+\lambda \frac{R^{2(n-1)}}{r^{n-1}(t)}v_k=0$$
Since
 $$\frac{R^{2(n-1)}}{r^{n-1}(t)}\geq R^{n-1}\geq r^{n-1}(t)$$
for all $t\geq t_k$, we can apply the theorem \ref{teorema9} to $u$ and $v_k$ and establish that
 $$\frac{u'(t)}{u(t)}\geq \frac{v'(t)}{v(t)},\;\; \;\;t\in(t_k,\tilde t_k).$$
 So, taking $\epsilon>0$ and integrating the inequality above from $t_k+\epsilon$ to $t$, we get
 $$\log \left(\frac{|u(t)|}{|u(t_k+\epsilon)|}\right) \geq \log\left(\frac{|v_{k}(t)|}{|v_{k}(t_k+\epsilon)|}\right)$$

 $$\frac{|u(t)|}{|v_{k}(t)|}\geq \frac {|u(t_k+\epsilon)|}{|v_k(t_k+\epsilon)|}$$
 Sending $\epsilon \rightarrow 0$ and using that $u'(t_k)=v'_{k}(t_k)\neq0$, the result to $v_k$ follows.
 The proof of the other inequality  follows the same ideas.

\begin{thm}
\label{extens}
Any solution $u$ of \,(\ref{sturm})\,on a interval $[t_0,t_0 + \delta]$ with initial value
$u(t_0)=x_0$, $u'(t_0)=x_1$ can be extended to $[t_0,\infty)$.
\end{thm}
Again, the proof is presented on \cite[theorem 2.2]{Manfredo}

\begin{thm}
\label{thm12}
For every $\lambda>0$ there is a function $u:[t_0,\infty)\mapsto\mathbb{R}$ for which $[r^{n-1}(t)u'(t)]'+\lambda\,r^{n-1}(t)u(t)=0$
 for every $t\geq t_0$.
\end{thm}
\textbf{Proof}

The ideia of the proof is to define a function which sends  each $T>t_0$ to the first eigenvalue of the  Sturm-Liouville problem below
 $$\left\{\begin{array}{l}
 (r^{n-1}(t)u')'+ r^{n-1}(t)\lambda_{T} u=0 \\ \\
 u(t_0)=u(T)=0

\end{array}\right.$$

  and conclude that this  function is onto $(0,\infty)$.
  Take $u$ as above and $v_s$ as a solution of

$$\left\{\begin{array}{l}
 (r^{n-1}(t)v')'+ r^{n-1}(t)\lambda_{T+s} v=0 \\ \\
 v(t_0)=v(T+s)=0

\end{array}\right.$$

Multiplying the first equation by $v_s$, the second one by $u$ and subtracting the results, we have
$$[r^{n-1}(t)(u'v_s-uv_s')]'+(\lambda_{T}-\lambda_{T+s})r^{n-1}(t)uv_s=0$$
Integrating from $t_0$ to $T$ we find
\begin{equation}
\label{dependcontinua}
 r^{n-1}(T)(u'v_s)(T)+ (\lambda_{T}-\lambda_{T+s})\int_{t_0}^{T}r^{n-1}(t)uv_s ds=0
\end{equation}
Now, notice that, by the continuous dependence of parameters, $$\lim_{s\to 0}\! v_s(T)=u(T)=0$$
So $$\lim_{s\to 0}(\lambda_{T}-\lambda_{T+s})=0$$ and then $\lambda_{T}$ is a continuous
function of $T$. From (\ref{dependcontinua}) we conclude that $\lambda_{T}$ is a decreasing function of $T$, since $u'(T)<0$ and $v_s(T)>0$.
In order to finish the proof, define
$$w(t)=\sin(\alpha\int_{t_0}^{t}\frac{ds}{r^{n-1}(s)})$$ where
$$\alpha=\frac{\pi}{\int_{t_0}^{T_{0}+T}\frac{ds}{r^{n-1}(s)}}$$
Then,
$$\left\{\begin{array}{l}
 (r^{n-1}(t)w')'+ \frac{\alpha^{2}}{r^{n-1}(t)}w(t)=0 \\ \\
 w(t_0)=w(T+T_{0})=0

\end{array}\right.$$

Using a completely analogous idea as above, it's immediate to check that
$$\int_{t_0}^{t_0+T}(\lambda_{T}r^{n-1}-\frac{\alpha^{2}}{r^{n-1}})uw=0$$
From the fact that $r$ is increasing, one concludes that
$$ (\lambda_{T}r^{n-1}(t_0)-\frac{\alpha^{2}}{r^{n-1}(T_{0})})\int_{t_0}^{t_0+T}uw\leq0$$
Consequently, $$0<\lambda_{T}\leq\frac{\alpha^{2}(T)}{r^{2(n-1)}(T_0)}$$
Since $\lim_{T\to\infty}\alpha(T)=0$,
 we have $$\lim_{T\to\infty}\lambda(T)=0$$
 A similar argument shows that $\lambda_{T}\rightarrow \infty$ when $T\rightarrow 0$.
 Then, $\lambda_{T}$ is onto $(0,\infty)$ and using the theorem \ref{extens} we reach the desired result.

\end{document}